\documentclass[reqno,12pt]{amsart}
\usepackage{amssymb}

\usepackage{amsmath}

\newtheorem{theorem}{Theorem}
\newtheorem{corollary}{Corollary}

\newtheorem{lemma}{Lemma}

\begin{document}

\title[Moderate deviations for the $L_1$-norm of kernel density estimators]
{Moderate deviations for the $L_1$-norm of kernel density estimators}

\author[A.Yu. Zaitsev]{Andrei Yu. Zaitsev}

\begin{abstract}The rate of normal approximation for the integral norm of kernel density estimators is investigated in the case of densities with power-type singularities. The quantities from the formulations of published results by the author are estimated. By assumption, the density tends to zero as a power-type function when the argument tends to infinity. Moreover, the density may have a finite number of power-type zeroes and of points with power-type tending to infinity. For such densities the size of zones of moderate deviations are found.\end{abstract}

\subjclass {Primary 60G50; secondary 60F05, 60F15, 60F17, 62G07}

\email{zaitsev@pdmi.ras.ru}
\address{St.~Petersburg Department of Steklov Mathematical Institute
\newline\indent
Fontanka 27,
\newline\indent
St.~Petersburg, 191023, Russia
\newline\indent
and
\newline\indent
St. Petersburg State University,
\newline\indent
Universitetskaya nab. 7/9,
\newline\indent
St. Petersburg, 199034, Russia}

\maketitle

\noindent Let $X,$ $X_{1},$ $X_{2},\ldots $ be a sequence of i.i.d.\ random
variables in $\mathbf{R}$ with density $f$. Let $\{h_{n}\}_{n\geq 1}$ be a
sequence of positive constants such that $h_{n}\rightarrow 0$ as $%
n\rightarrow \infty $. The classical kernel estimator is defined as
\[
f_{n}(x)\stackrel{\mathrm{def}}{=}\frac{1}{nh_{n}}\sum_{i=1}^{n}K\left(
\frac{x-X_{i}}{h_{n}}\right) ,\quad \mbox{for}\;x\in \mathbf{R},
\]
where $K$ is a kernel satisfying $K(u)=0$, for $|u|>1/2$, $$\|K\|_{\infty
}=\sup_{u\in \mathbf{R}}\left| K(u)\right| =\kappa <\infty $$ and $\int_{%
\mathbf{R}}K(u)\,du=1$. Let $\|\cdot \|$ denote the $L_{1}(\mathbf{R})$%
-norm. Write $\|K^{2}\|=\int_{\mathbf{R}}K^{2}(u)\,du$.

Devroye and Gy\"{o}rfi \cite{dg} posed the question about the asymptotic
distribution of $\|f_{n}-f\|.$

Hall \cite{ha}, M. Cs\"{o}rg\H{o} and Horv\'{a}th \cite{ch} and Horv\'{a}th \cite{h}
proved the Central Limit Theorem (CLT) for $\|f_{n}-f\|_{p},$ the $L_{p}$-norm distance, $p\geq 1$, under some additional regularity restrictions on
the density $f$. Horv\'{a}th \cite{h} introduced a Poissonization technique
into the study of CLTs for $\|f_{n}-f\|_{p}$. Beirlant and Mason \cite{bm}
proposed a general method for proving the asymptotic normality of the $L_{p}$%
-norm of empirical functionals. Mason (see Theorem 8.9 in Eggermont and
LaRiccia \cite{er}) has applied their method to the special case of the $%
L_{1}$-norm of the kernel density estimator and proved Theorem \ref{t1.1}
below. Gin\'{e}, Mason and Zaitsev \cite{gmz} extended the CLT result of
Theorem \ref{t1.1} to processes indexed by kernels~$K$. One should mention the papers by 
Louani \cite{L} and Gao \cite{G} in which the logarithmic asymptotics of the probabilities of moderate deviations is investigated. Similar questions are considered in \cite{df} and \cite{LWX}.

\medskip

 \begin{theorem} [Mason (see Theorem 8.9 in Eggermont and LaRiccia
\cite{er}]\label{t1.1} For any Lebesgue density $f$ and for any sequence of positive
constants $\left\{ h_{n}\right\} _{n\geq 1}$ satisfying $h_{n}\rightarrow 0$
and $nh_{n}^{2}\rightarrow \infty $, as $n\rightarrow \infty$, we have
\[
\frac{\Vert f_{n}-\mathbf{E\,}f_{n}\Vert -\mathbf{\,E\,}\Vert f_{n}-\mathbf{%
\,E\,}f_{n}\Vert }{\sqrt{\mbox{\rm Var}(\Vert f_{n}-\mathbf{\,E\,}f_{n}\Vert
)}}\rightarrow _{d}Z
\]
and
\[
\lim_{n\rightarrow \infty }n\,\mbox{{\rm Var}}(\Vert f_{n}-\mathbf{E\,}
f_{n}\Vert )=\sigma ^{2},
\]
where
\[
\sigma ^{2}=\sigma ^{2}(K)\stackrel{\mathrm{def}}{=}\|K^{2}\|\int_{-1}^{1}%
\mbox{{\rm cov}}\left( \left\vert \sqrt{1-\rho ^{2}(t)}\,Z+\rho
(t)\,Y\right\vert ,\left\vert Y\right\vert \right) dt,
\]
$Z\ $and $Y$ are independent standard normal random variables and, for any $
t\in {\mathbf{R}},$
\[
\rho (t)=\rho (t,K)\stackrel{\mathrm{def}}{=}\frac{\int_{\mathbf{R}%
}K(u)\,K(u+t)\,du}{\|K^{2}\|}.
\] \end{theorem}
\medskip

It is well-known that $\rho (t)$ is the characteristic function of a
symmetric probability distribution (see, e.g., Lukacs \cite{luc}, Theorem
4.2.4). Clearly, $\rho (t)$ is a continuous function of $t$, $\left\vert
\rho (t)\right\vert \leq 1$, $\rho (0)=1$ and $\rho (t)=0$ for $|t|\geq 1$.
The variance $\sigma ^{2}$ has an alternate representation.
\[
\mbox{{\rm cov}}\left( \left\vert \sqrt{1-\rho ^{2}(t)}\,Z+\rho
(t)\,Y\right\vert ,\left\vert Y\right\vert \right) =\varphi \left( \rho
(t)\right) ,
\]
where
\begin{equation}
\varphi (\rho )\stackrel{\mathrm{def}}{=}\frac{2}{\pi }\left( \rho \arcsin
\rho +\sqrt{1-\rho ^{2}}-1\right) ,\quad \rho \in \left[ -1,1\right] .
\label{phiofrho}
\end{equation}

Theorem \ref{t1.1} shows that $\Vert f_{n}-\mathbf{E\,}f_{n}\mathbf{\Vert }$
is asymptotically normal under no assumptions at all on the density $f$.
Centering by $\mathbf{E\,}f_{n}$ is more natural from a probabilistic point
of view. The estimation of $\Vert f-\mathbf{E\,}f_{n}\mathbf{\Vert }$ (if
needed) is a purely analytic problem. The main results of Zaitsev \cite{zai}
(Theorems \ref{t1.2}, \ref{t1.3} and \ref{t1.4}) provide estimates of the
rate of strong approximation and bounds for probabilities of moderate
deviations in the CLT of Theorem \ref{t1.1}. These results are general and
provides estimates for arbitrary densities. Therefore, these results have
complicated formulations which are expressed in terms of some sets $E_{n}$.
Roughly speaking, $\mathbf{\,}f_{n}$ is close to $f$ on the set $E_{n}$. To
obtain the rate of approximation in terms of $h_{n}$ one needs additional
considerations. In Zaitsev \cite{zai}, one can find a number of concrete
examples of calculating such a rate.

The aim of this paper is to obtain the rates of approximation for densities
with power-type singularities. We provide the bound for the quantities
involved in the formulations of the results of Zaitsev \cite{zai} assuming
that the density  $f(x)$ tends to zero as a power-type function as $%
x\rightarrow \pm \infty $, and has a finite number of power-type zeroes and
of the points with power-type tending to infinity. For such density we shal
find the size for zones of moderate deviations.

Le us formulate the results of Zaitsev \cite{zai}. Set, for any Borel sets $%
B,E$,
\[
R_{n}(B,E)\stackrel{\mathrm{def}}{=}\int_{B}\left( \int_{-1}^{1}\left|
g_{n}(x,t,E)-g(x,t,E)\right| \,dt\right) \,dx,
\]
where
\[
g(x,t,E)\stackrel{\mathrm{def}}{=}\mathbf{1}_{E}(x)\,\mbox{{\rm cov}}\left(
\left| \sqrt{1-\rho ^{2}(t)}\,Z+\rho (t)\,Y\right| ,\left| Y\right| \right)
\,f(x),
\]
\[
g_{n}(x,t,E)\stackrel{\mathrm{def}}{=}\mathbf{1}_{E}(x)\mathbf{1}%
_{E}(x+th_{n})\,\mathbb{C}_{n}\left( x,x+th_{n}\right) \,\sqrt{%
f(x)\,f(x+th_{n})},
\]
\[
\mathbb{C}_{n}\left( x,y\right) \stackrel{\mathrm{def}}{=}\mbox{{\rm cov}}%
\left( \left| \sqrt{1-\rho _{n,x,y}^{2}}Z+\rho _{n,x,y}\,Y\right|
,|Y|\right) ,
\]
$Z$\ and $Y$\ are independent standard normal random variables and
\[
\rho _{n,x,y}\stackrel{\mathrm{def}}{=}\frac{\mathbf{E\,}\left[ K\left(
\frac{x-X}{h_{n}}\right) \,K\left( \frac{y-X}{h_{n}}\right) \right] }{\sqrt{%
\mathbf{E\,}K^{2}\left( \frac{x-X}{h_{n}}\right) \mathbf{\,E\,}K^{2}\left(
\frac{y-X}{h_{n}}\right) }}.
\]

The following Lemma \ref{L1} is crucial for the formulation of the main results,
Theorems \ref{t1.2}--\ref{t1.4} below.\medskip

\begin{lemma} \label{L1} Whenever $h_{n}\rightarrow 0$\ and $%
nh_{n}^{2}\rightarrow \infty $, as $n\rightarrow \infty $,\ there exist
sequences of Borel sets
\[
E_{1}\subset E_{2}\subset \cdots \subset E_{n}\subset \cdots
\]
and constants $\left\{ \beta _{n}\right\} _{n=1}^{\infty }$\ and $\left\{
D_{n}\right\} _{n=1}^{\infty }$ such that the relations
\begin{equation}
\phi _{n}\stackrel{\mathrm{def}}{=}\int_{\mathbf{R}\backslash
E_{n}}f(x)\,dx\rightarrow 0,\quad \mbox{{as} }n\rightarrow \infty ,
\label{psin0}
\end{equation}
\begin{equation}
0<\beta _{n}\stackrel{\mathrm{def}}{=}\inf_{y\in E_{n}}f(y)\leq f(x)\leq
D_{n}\stackrel{\mathrm{def}}{=}\sup_{y\in E_{n}}f(y)<\infty ,\quad %
\mbox{{for} }x\in E_{n},  \label{bdn}
\end{equation}
and
\begin{equation}
\varepsilon _{n}\stackrel{\mathrm{def}}{=}\sup_{H\in \mathcal{H}%
_{0}}\;\sup_{x\in E_{n}}\left| h_{n}^{-1}\int_{\mathbf{R}}f(z)\,H\left(
\frac{x-z}{h_{n}}\right) \,dz-f(x)\int_{\mathbf{R}}H(z)\,dz\right|
\rightarrow 0,  \label{CC}
\end{equation}
as $n\rightarrow \infty $, are valid, where
\begin{equation}
\mathcal{H}_{0}\stackrel{\mathrm{def}}{=}\left\{ K,K^{2},\left| K\right|
^{3},\mathbf{1}\{x:|x|\leq 1/2\}\right\} .  \label{H0}
\end{equation}
Moreover,
\begin{eqnarray}
\frac{D_{n}^{1/2}}{\beta _{n}^{1/2}}\left( \frac{1}{\left( \beta
_{n}\,nh_{n}\right) ^{1/5}}+\frac{\varepsilon _{n}}{\beta _{n}}\right)
&+&R_{n}(E_{n},E_{n})+\frac{\lambda (E_{n})}{\sqrt{nh_{n}^{2}}}+D_{n}\,h_{n}
\nonumber \\
&+&\frac{D_{n}^{3}\,P_{n}}{\beta _{n}^{3}}+\mathbb{N}_{n}\sqrt{h_{n}}%
\rightarrow 0,  \label{tendzero}
\end{eqnarray}
as\ $n\rightarrow \infty $, where $\lambda (\,\cdot \,)$\ means the Lebesgue
measure,
\[
\mathbb{N}_{n}\stackrel{\mathrm{def}}{=}\int_{E_{n}}f^{3/2}(x)\,dx\ \quad
\text{and}\ \quad P_{n}\stackrel{\mathrm{def}}{=}\max_{x\in \mathbf{R}}%
\mathbf{P}\left\{ [x,x+2\,h_{n}]\right\} .
\]\end{lemma}
\medskip

\begin{theorem} \label{t1.2} There exists an absolute constant $A$\ such
that, whenever $h_{n}\rightarrow 0$\ and $nh_{n}^{2}\rightarrow \infty $, as
$n\rightarrow \infty $, for any sequence of Borel sets $E_{1},E_{2},\ldots
,E_{n},\ldots $\ satisfying $(\ref{psin0})$--$(\ref{tendzero})$, there exists an
$n_{0}\in N$ such that, for any fixed $x>0$\ and for sufficiently large
fixed $n\geq n_{0}$, one can construct on a probability space a sequence of
i.i.d.\ random variables $X_{1},X_{2},\ldots $ and a standard normal random
variable $Z$ such that
\begin{eqnarray*}
&&\hspace{1cm}\mathbf{P}\left\{ \left| \sqrt{n}\mathbf{\,}\Vert f_{n}-%
\mathbf{E\,}f_{n}\Vert -\sqrt{n}\mathbf{\,E\,}\Vert f_{n}-\mathbf{E\,}%
f_{n}\Vert -\sigma \,Z\right| \geq y_{n}+z+x\right\} \\
&\leq &A\,\Big(\exp \left\{ -A^{-1}\,\sigma ^{-1}x/\tau _{n}\right\} \\
&&+\exp \left\{ -A^{-1}\,\kappa ^{-1}\,\Omega _{n}^{-1/2}z\,\log ^{*}\log
^{*}(z/A\,\kappa \,\Omega _{n}^{1/2})\right\} \\
&+&\mathbf{P}\left\{ \left| \delta _{n}Z\right| \geq z/2\right\} \Big),\quad %
\mbox{{for any} }z>0,
\end{eqnarray*}
where $\log ^{*}b=\max \left\{ 1,\log b\right\} $,
\[
\tau _{n}\stackrel{\mathrm{def}}{=}A\,\Psi _{n}^{3/2}\,\left( P_{n}+\psi
_{n}\right) ^{1/2}\rightarrow 0,\mathbf{\quad }\mbox{{as} }n\rightarrow
\infty ,
\]
\[
y_{n}\stackrel{\mathrm{def}}{=}\frac{A\,\lambda (E_{n})\,\left\|
K^{3}\right\| }{\left\| K^{2}\right\| \sqrt{nh_{n}^{2}}}+\frac{A\,\mathbb{N}%
_{n}\sqrt{h_{n}}}{\sqrt{\left\| K^{2}\right\| }}\rightarrow 0,\mathbf{\quad }%
\mbox{{as} }n\rightarrow \infty ,
\]
\begin{eqnarray*}
\delta _{n} &\stackrel{\mathrm{def}}{=}&\frac{A\,\|K^{2}\|}{\sigma \,h_{n}}%
\,\left( \mathbb{L}_{n}+\frac{\varepsilon _{n}\,\mathbb{M}_{n}}{\left\|
K^{2}\right\| }\right) \\
&+&A\,\kappa \,\Omega _{n}^{1/2}+\frac{A}{\sigma }\left( \frac{\left\|
K^{3}\right\| \,\lambda (E_{n})}{\left\| K^{2}\right\| \sqrt{nh_{n}^{2}}}%
\right) ^{2}\rightarrow 0,\mathbf{\quad }\mbox{{as} }n\rightarrow \infty ,
\end{eqnarray*}
\[
\mathbb{L}_{n}\stackrel{\mathrm{def}}{=}\int_{E_{n}}\int_{E_{n}}\mathbf{1}%
\{|x-y|\leq h_{n}\}\,\sqrt{f(x)\,f(y)}\,\,\mathbb{K}_{n}(x,y)\,dx\,dy,
\]
\[
\mathbb{K}_{n}(x,y)\stackrel{\mathrm{def}}{=}\min \left\{ 1-\rho
_{n,x,y}^{2},\frac{\left\| K^{3}\right\| }{\left( 1-\rho _{n,x,y}^{2}\right)
^{3/2}\left\| K^{2}\right\| ^{3/2}\sqrt{n\,h_{n}\,f(x)}}\right\},
\]
\[
\mathbb{M}_{n}\stackrel{\mathrm{def}}{=}\int_{E_{n}}\int_{E_{n}}\mathbf{1}%
\{|x-y|\leq h_{n}\}\,f^{1/2}(x)\,f^{-1/2}(y)\,dx\,dy,
\]
\[
\Omega _{n}\stackrel{\mathrm{def}}{=}\alpha _{n}+2\,P_{n}+2\,\phi _{n}+\frac{%
4\,\|K^{2}\|\,R_{n}(E_{n},E_{n})}{\sigma ^{2}}+L(n,\mathbf{R)}\rightarrow 0,%
\mathbf{\quad }\mbox{{as} }n\rightarrow \infty \mathbf{,}
\]
\[
\alpha _{n}\stackrel{\mathrm{def}}{=}\frac{1296}{5}\tau _{n}^{2}\log \frac{1%
}{\tau _{n}},\mathbf{\quad }\Psi _{n}\stackrel{\mathrm{def}}{=}\left\|
K^{2}\right\| \,D_{n}\,\beta _{n}^{-1}\,\kappa ^{2}\,\sigma ^{-4},
\]
\[
\psi _{n}\stackrel{\mathrm{def}}{=}256\,\kappa ^{2}\,\sigma ^{-2}\,\min
\left\{ P_{n},D_{n}\,h_{n}\right\} ,
\]
\[
L(n,\mathbf{R})\stackrel{\mathrm{def}}{=}\int_{\mathbf{R}}\left| h_{n}^{-1}\,%
\mathbf{P}\{X\in [x-h_{n}/2,x+h_{n}/2]\}-f(x)\right| \,dx.
\]

Denote by $F\{\,\cdot \,\}$ and $\Phi \{\,\cdot \,\}$ the distributions of
the random variables $$\,\sqrt{n}\,\left( \Vert f_{n}-\mathbf{E}\,f_{n}\Vert -
\mathbf{E\,}\Vert f_{n}-\mathbf{E\,}f_{n}\Vert \right) /\sigma $$ and $\,Z$,
respectively. The Prokhorov distance is defined by $$\pi (F,\Phi )=\inf
\left\{ \varepsilon :\pi (F,\Phi ,\varepsilon )\leq \varepsilon \right\}, $$
where
\[
\pi (F,\Phi ,\varepsilon )=\sup_{X}\max \left\{ F\{X\}-\Phi \{X^{\varepsilon
}\}\mbox{, }\Phi \{X\}-F\{X^{\varepsilon }\}\right\} ,\quad \varepsilon >0,
\]
and $X^{\varepsilon }$ is the $\varepsilon $-neighborhood of the Borel set $%
X $.\end{theorem} \medskip

\begin{corollary} \label{c1.1} There exists an absolute constant $A$\ such
that, whenever $h_{n}\rightarrow 0$\ and $nh_{n}^{2}\rightarrow \infty $, as
$n\rightarrow \infty $, for any sequence of Borel sets $E_{1},E_{2},\ldots
,E_{n},\ldots $\ satisfying $(\ref{psin0})$--$(\ref{tendzero})$, there exists an
$n_{0}\in N$ such that, for sufficiently large fixed $n\geq n_{0}$ and for
any $\varepsilon >0$,
\begin{eqnarray*}
&&\pi (F,\Phi ,2\,\varepsilon +y_{n}/\sigma ) \\
&\leq &A\,\left( \exp \left\{ -A^{-1}\,\kappa ^{-1}\,\Omega
_{n}^{-1/2}\sigma \,\varepsilon \,\log ^{*}\log ^{*}(\sigma \,\varepsilon
/A\,\kappa \,\Omega _{n}^{1/2})\right\} \right.
\end{eqnarray*}
\[
\left. +\;\exp \left\{ -A^{-1}\,\varepsilon /\tau _{n}\right\} +\mathbf{P}%
\left\{ \left| \delta _{n}Z\right| \geq \sigma \,\varepsilon /2\right\}
\right)
\]
and
\[
\mathbb{\pi }(F,\Phi )\leq y_{n}/\sigma +A\,\tau _{n}\,\log ^{*}\left(
1/\tau _{n}\right) +A\,\delta _{n}\,\sigma ^{-1}\,\sqrt{\log ^{*}\left(
\sigma /\delta _{n}\right) }
\]
\[
+\,A\,\kappa \,\Omega _{n}^{1/2}\sigma ^{-1}\,\log ^{*}\left( \sigma /\kappa
\,\Omega _{n}^{1/2}\right) /\log ^{*}\log ^{*}(\sigma /\kappa \,\Omega
_{n}^{1/2}).
\]
\end{corollary}\medskip

\begin{theorem} \label{t1.3} There exists an absolute constant $A$\ such
that, whenever $h_{n}\rightarrow 0$\ and $nh_{n}^{2}\rightarrow \infty $, as
$n\rightarrow \infty $, for any sequence of Borel sets $E_{1},E_{2},\ldots
,E_{n},\ldots $\ satisfying $(\ref{psin0})$--$(\ref{tendzero})$,\ there exists
an $n_{0}\in N$ such that,\ for sufficiently large fixed $n\geq n_{0}$\ and
for any fixed $b$\ satisfying $\tau _{n}\leq A^{-1}b$, $b\leq 1$, one can
construct on a probability space a sequence of i.i.d.\ random variables $%
X_{1},X_{2},\ldots $ and a standard normal random variable $Z$\ such that
\[
\hspace{-1cm}\mathbf{P}\left\{ \left| \sqrt{n}\mathbf{\,}\Vert f_{n}-\mathbf{%
E\,}f_{n}\Vert -\sqrt{n}\mathbf{\,E\,}\Vert f_{n}-\mathbf{E\,}f_{n}\Vert
-\sigma \,Z\right| \right.
\]
\[
\hspace{1cm}\geq \left. A\,\sigma \,\exp \{-b^{2}/72\left( \tau _{n}\right)
^{2}\}+y_{n}+z+x\right\}
\]
\begin{eqnarray*}
&\leq &A\,\Big(\exp \left\{ -A^{-1}\,\sigma ^{-1}x/\tau _{n}\right\} \\
&&+\exp \{-A^{-1}\,\kappa ^{-1}\,\Omega _{n}^{-1/2}z\,\log ^{*}\log
^{*}(z/A\,\kappa \,\Omega _{n}^{1/2})\}
\end{eqnarray*}
\[
{}+\mathbf{P}\left\{ b\,\left| Z\right| >A^{-1}\,\sigma ^{-1}x\right\} +%
\mathbf{P}\left\{ \left| \delta _{n}Z\right| \geq z/2\right\} \Big),\quad %
\mbox{{for any} }x,z>0.
\]\end{theorem}
\medskip

In the formulations of Theorems \ref{t1.2} and \ref{t1.3} and Corollary \ref{c1.1}, the numbers $n_{0}$
depend on $\{h_{n}\}_{n\geq 1}$, $\{E_{n}\}_{n\geq 1}$, $f$ and $K$.

Denote now by $F(\,\cdot \,)$\ and $\Phi (\,\cdot \,)$\ the distribution
functions of random variables $\sqrt{n}\,\left( \Vert f_{n}-\mathbf{E\,}%
f_{n}\Vert -\mathbf{E\,}\Vert f_{n}-\mathbf{E\,}f_{n}\Vert \right) /\sigma $
and $\,Z$,\ respectively. For example, $\Phi (x)=\Phi \left\{ (-\infty
,x]\right\} $. The following statement about moderate deviations follows
from Theorem~\ref{t1.2}.\medskip

\begin{theorem} \label{t1.4} Under the conditions of Theorem $\ref{t1.2}$, we have
\[
\frac{F(-x)}{\Phi (-x)}\rightarrow 1\mathit{\ }\quad \text{and}\mathit{\ }%
\quad \frac{1-F(x)}{1-\Phi (x)}\rightarrow 1,\mathbf{\quad }\mbox{{as} }%
n\rightarrow \infty ,
\]
if
\[
0<x=x_{n}=o\left( \min \left\{ \tau _{n}^{-1/3},\mathbf{\,\,\,}\Omega
_{n}^{-1/6}\left( \log ^{*}\log ^{*}(1/\,\Omega _{n})\right) ^{1/3},\mathbf{%
\,\,\,}y_{n}^{-1},\mathbf{\,\,\,}\delta _{n}^{-1/2}\right\} \right) .
\]
\end{theorem}

The choice of sets $E_{n}$, which are involved in the formulations of our
results, is not unique. Lemma 1 ensures that, for \textit{any} density $f$,
there exist sets $E_{n}$ such that the quantities $\tau _{n},y_{n},\Omega
_{n}$ and $\delta _{n}$ tend to zero. The optimization of the choice of $%
E_{n}$ is a separate problem. However, for sufficiently regular densities $f$%
, it is not difficult to choose $E_{n}$ so that the rate of approximation is
good enough, see the examples below.

Below we consider a number of examples with calculating the concrete order
of dependence of our estimates with respect to $h_{n}$. For this we shall
need some simplifying restriction on $K$ and $f$.

We shall assume that the kernel $K$ is such that
\begin{equation}
\rho (t)\leq 1-\left( C_{0}|t|\right) ^{d},\quad \mbox{for }|t|\leq 1,
\label{rh2}
\end{equation}
with some $0<C_{0}\leq 1$, $1\leq d\leq 2.$\ Using that $\rho (t)$ is a
characteristic function of a probability distribution, it is easy to show
that each kernel considered here satisfies this condition with $d=2$ and
some $C_{0}$. But for some standard kernels such a condition is valid for
smaller $d$. For example, for $K(x)=\mathbf{1}\{x:|x|\leq 1/2\}$, inequality
(\ref{rh2}) is valid with $C_{0}=1$ and $d=1$. Note that there exist kernels
satisfying (\ref{rh2}) with $0<d<1$ but, for the most of natural kernels, (%
\ref{rh2}) is valid with some $d$ such that $1\leq d\leq 2.$

Denote
\begin{equation}
q_{n}(x)\stackrel{\mathrm{def}}{=}\sup_{|z-x|\leq 3h_{n}/2}\left\vert
f(z)-f(x)\right\vert .  \label{7a}
\end{equation}
In our examples, the choice of $E_{n}$ will be such that
\begin{equation}
Q_{n}\stackrel{\mathrm{def}}{=}\sup_{x\in E_{n}}\frac{q_{n}(x)}{f(x)}%
\rightarrow 0,\quad \mbox{
as }n\rightarrow \infty ,  \label{rh4}
\end{equation}
and
\begin{equation}
\sup_{x\in E_{n}}\Delta _{n}(x)\rightarrow 0,\quad \mbox{
as }n\rightarrow \infty ,  \label{rh7}
\end{equation}
where
\begin{equation}
\Delta _{n}(x)\stackrel{\mathrm{def}}{=}\frac{\left\Vert K^{3}\right\Vert }{%
\left\Vert K^{2}\right\Vert ^{3/2}\sqrt{n\,h_{n}\,f(x)}}.  \label{DELTA}
\end{equation}
Note that the sets $E_{n}$ have to satisfy the conditions (\ref{psin0}), (%
\ref{CC}) and (\ref{tendzero}) if Theorems \ref{t1.2} and \ref{t1.3} and Corollary \ref{c1.1} are to be
applied.

Below we shall estimate the quantities appearing in our theorems using
simplifying conditions (\ref{rh2}) and (\ref{rh4}). Denote
\begin{equation}
G_{n}(x)\stackrel{\mathrm{def}}{=}\sup_{|z-x|\leq h_{n}}f(z)  \label{rh3}
\end{equation}
and
\[
B_{n}(x)\stackrel{\mathrm{def}}{=}\inf_{|z-x|\leq h_{n}}f(z).
\]
By (\ref{rh4}), we have
\begin{equation}
G_{n}(x)\leq 2\,f(x)\quad \mbox{for }x\in E,  \label{rh6}
\end{equation}
for sufficiently large $n$.

It is easy to see that
\begin{equation}
\mathbf{E\,}h_{n}^{-1}\,H\left( \frac{x-X}{h_{n}}\right) =h_{n}^{-1}\int_{%
\mathbf{R}}f(z)\,H\left( \frac{x-z}{h_{n}}\right) \,dz=\int_{\mathbf{R}%
}f(x-h_{n}y)\,H(y)\,dy  \label{Sy}
\end{equation}
and if $H(x)=0$, for $|x|>3/2$, then
\begin{equation}
\left\vert h_{n}^{-1}\int_{\mathbf{R}}f(z)\,H\left( \frac{x-z}{h_{n}}\right)
\,dz-f(x)\int_{\mathbf{R}}H(y)\,dy\right\vert \leq q_{n}(x)\,\int_{\mathbf{R}%
}\left\vert H(y)\right\vert \,dy,  \label{SV}
\end{equation}
This bound allows us to estimate $\varepsilon _{n}$ and $L(n,\mathbf{R})$.

Recall that
\[
g(x,t,E)=\mathbf{1}_{E}(x)\,\varphi \left( \rho (t)\right) \,f(x),
\]
\[
g_{n}(x,t,E)=\mathbf{1}_{E}(x)\mathbf{1}_{E}(x+th_{n})\,\varphi \left( \rho
_{n,x,y}\right) \,\sqrt{f(x)\,f(x+th_{n})},
\]
where $y=x+th_{n}$. In our treatment of the examples, we shall use the fact
that the function $\varphi (\rho )$ in (\ref{phiofrho}) satisfies $0\leq
\varphi (\rho )\leq 1$ and the Lipschitz condition $\left\vert \varphi (\rho
_{1})-\varphi (\rho _{2})\right\vert \leq \left\vert \rho _{1}-\rho
_{2}\right\vert $. The most complicated estimation is that of
\begin{eqnarray}
R_{n}(E_{n},E_{n}) &=&\int_{E_{n}}\left( \int_{-1}^{1}\left\vert
g_{n}(x,t,E_{n})-g(x,t,E_{n})\right\vert \,dt\right) \,dx  \nonumber \\
&=&\int_{E_{n}^{\ast }}+\int_{E_{n}\backslash E_{n}^{\ast }}=I_{1}+I_{2},
\label{R1}
\end{eqnarray}
where $E_{n}^{\ast }=\left\{ x\in E_{n}:x+th_{n}\in E_{n}\mbox{ for all }%
t:|t|\leq 1\right\} $. The set $E_{n}\backslash E_{n}^{\ast }$ contains
points belonging to $E_{n}$ and to the $h_{n}$-neighborhood of the boundary
of $E_{n}$. Integrating over $E_{n}\backslash E_{n}^{\ast }$, we use the
trivial inequality
\begin{eqnarray*}
&&\left\vert g_{n}(x,t,E_{n})-g(x,t,E_{n})\right\vert \\
&\leq &\max \left\{ \mathbf{1}_{E_{n}}(x)\,\,f(x),\;\mathbf{1}_{E_{n}}(x)%
\mathbf{1}_{E_{n}}(x+th_{n})\,\sqrt{f(x)\,f(x+th_{n})}\right\} ,
\end{eqnarray*}
and
\begin{equation}
I_{2}\leq \lambda (E_{n}\backslash E_{n}^{\ast })\max f(y),  \label{R2}
\end{equation}
where the maximum is taken over points of the form $y=x+th_{n}\in E_{n}$
such that $x\in E_{n}\backslash E_{n}^{\ast }$ and $|t|\leq 1$. For $x\in
E_{n}^{\ast }$, $|t|\leq 1$, we have
\begin{eqnarray*}
&&\left\vert g_{n}(x,t,E_{n})-g(x,t,E_{n})\right\vert \\
&\leq &\left\vert \varphi \left( \rho (t)\right) -\varphi \left( \rho
_{n,x,y}\right) \right\vert \,f(x)+\varphi \left( \rho _{n,x,y}\right)
\left\vert f(x)-\sqrt{f(x)\,f(y)}\right\vert \\
&\leq &\left\vert \rho (t)-\rho _{n,x,y}\right\vert \,f(x)+\left\vert f(x)-%
\sqrt{f(x)\,f(y)}\right\vert ,
\end{eqnarray*}
where $y=x+th_{n}$. Furthemore,
\[
\rho _{n,x,y}=\frac{\mathbf{E\,}\left[ K\left( \frac{x-X}{h_{n}}\right)
\,K\left( \frac{x-X}{h_{n}}+t\right) \right] }{\sqrt{\mathbf{E\,}K^{2}\left(
\frac{x-X}{h_{n}}\right) \mathbf{\,E\,}K^{2}\left( \frac{x-X}{h_{n}}%
+t\right) }}.
\]
Using (\ref{rh4}) and applying (\ref{SV}) for $H(x)=K^{2}(x)$, $%
H(x)=K^{2}(x+t)$ and $H(x)=K(x)\,K(x+t)$, we obtain
\begin{equation}
\left\vert \rho (t)-\rho _{n,x,y}\right\vert \,f(x)\leq 3\,\sup_{|z-x|\leq
3h_{n}/2}\left\vert f(z)-f(x)\right\vert =3\,q_{n}(x),  \label{rh}
\end{equation}
for sufficiently large $n$. Moreover,
\[
\left\vert f(x)-\sqrt{f(x)\,f(y)}\right\vert \leq \sup_{|z-x|\leq
h_{n}}\left\vert f(z)-f(x)\right\vert \leq q_{n}(x),
\]
Thus,
\begin{equation}
I_{1}\leq 4\,\int_{E_{n}^{\ast }}q_{n}(x)\,dx,  \label{R3}
\end{equation}
for sufficiently large $n$.

By definition,
\begin{equation}
\mathbb{K}_{n}(x,y)=\min \left\{ 1-\rho _{n,x,y}^{2},\frac{\Delta }{\left(
1-\rho _{n,x,y}^{2}\right) ^{3/2}}\right\} ,  \label{KN}
\end{equation}
where $\Delta =\Delta _{n}(x)$. When estimating
\[
\mathbb{L}_{n}=h_{n}\int_{E_{n}}\left( \int_{-1}^{1}\mathbf{1}\{x+th_{n}\in
E_{n}\}\,\sqrt{f(x)\,f(x+th_{n})}\,\,\mathbb{K}_{n}(x,x+th_{n})\,dt\right)
\,dx
\]
in the examples, we will use the fact that, by (\ref{bdn}) and (\ref{KN}),
we have
\begin{equation}
\mathbb{K}_{n}(x,y)\leq \Delta ^{2/5}.  \label{KN0}
\end{equation}

Let us estimate the integral
\[
J=\int_{-1}^{1}\mathbf{1}\{x+th_{n}\in E_{n}\}\,\sqrt{f(x)\,f(x+th_{n})}\,\,%
\mathbb{K}_{n}(x,x+th_{n})\,dt.
\]
By (\ref{rh3}), is may be bounded from above by
\begin{eqnarray*}
G_{n}(x)\int_{-1}^{1}\mathbf{1}\{x+th_{n} &\in &E_{n}\}\,\mathbb{K}%
_{n}(x,x+th_{n})\,dt=G_{n}(x)\left( \int_{|t|\leq \delta }+\int_{\delta \leq
|t|\leq 1}\right) \\
&=&G_{n}(x)\left( J_{1}+J_{2}\right) ,
\end{eqnarray*}
where $\delta =\left( 6\,q_{n}(x)/f(x)+\Delta ^{2/5}\right) ^{1/d}/C_{0}$.
Using (\ref{rh2}), (\ref{rh}) and (\ref{KN0}), we see that if $2/3<d\leq 2$,
then
\begin{multline}
J_{1}+J_{2} \leq 2\,\delta\, \Delta ^{2/5}+\int_{\delta \leq |t|\leq 1}\frac{%
\Delta \,dt}{\left( \big( C_{0}|t|\right) ^{d}/2\big) ^{3/2}} \\
\leq 2\,\delta\, \Delta ^{2/5}+c(d)\,\Delta ^{\left( 2+2d\right)
/5d}/C_{0}\leq c(d)\,\Delta ^{2/5}\left( q_{n}(x)/f(x)+\Delta ^{2/5}\right)
^{1/d}/C_{0},
\end{multline}
for sufficiently large $n$. Hence, taking into account (\ref{rh6}), we see
that
\begin{equation}
\mathbb{L}_{n}\leq \frac{c(d)\,h_{n}}{C_{0}}\int_{E_{n}}\Delta ^{2/5}\left(
q_{n}(x)/f(x)+\Delta ^{2/5}\right) ^{1/d}f(x)\,dx,  \label{Ln}
\end{equation}
for sufficiently large $n$.

Considering the examples we assume always that the kernel $K$ satisfies (\ref
{rh2}) and we shall verify the validity of (\ref{psin0}), (\ref{CC}), (\ref
{tendzero}), (\ref{rh4}) and (\ref{rh7}). We shall denote $l_{n}=\log \log
\frac{1}{h_{n}}$.\medskip

\noindent {\bf Example 1.} Consider the density $f$ satisfying $f(x)=0$, for $%
x\notin \lbrack 0,1]$,
\[
C_{1}\left\vert x\right\vert ^{\lambda }\leq f(x)\leq C_{2}\left\vert
x\right\vert ^{\mu },\quad \mbox{for }0<x\leq 1,
\]
and the Lipschitz type condition
\begin{equation}
\left\vert f(x)-f(y)\right\vert \leq C_{3}\left\vert x-y\right\vert ^{\gamma
},\quad \quad \mbox{for }x,y\in \lbrack 0,1],  \label{LIP}
\end{equation}
where $C_{1}$, $C_{2}$, $C_{3}$, $0\leq \mu \leq \lambda <\infty $, and $%
0<\gamma \leq 1$, are constants. Choose
\[
E_{n}=\left[ u_{n},1-2h_{n}\right] ,
\]
where
\[
u_{n}=\max \left\{ \left( h_{n}\log \frac{1}{h_{n}}\right) ^{\alpha
},\varkappa _{n}\,h_{n}^{2\gamma /3\lambda },\varkappa _{n}\,\left(
nh_{n}\right) ^{-2/7\lambda }\right\} ,\quad \alpha =\frac{1}{1+3\lambda
+\mu },
\]
and $\varkappa _{n}\rightarrow \infty $ as $n\rightarrow \infty $ is an
arbitrary sequence of positive numbers which tend to infinity slower than a
fixed function of $h_{n}^{-1}$. Without loss of generality we assume $%
h_{n}\leq 1/8$. Then we estimate $\phi _{n}=O\left( u_{n}^{1+\mu }\right) $,
$\beta _{n}^{-1}=O\left( u_{n}^{-\lambda }\right) $, $D_{n}=O\left( 1\right)
$. Moreover, by (\ref{CC}), (\ref{7a}), (\ref{SV}) and (\ref{LIP}), we have $%
\sup_{x\in E_{n}}q_{n}(x)=O\left( h_{n}^{\gamma }\right) $, $\varepsilon
_{n}=O\left( h_{n}^{\gamma }\right) $, $Q_{n}=O\left( h_{n}^{\gamma
}u_{n}^{-\lambda }\right) $. The relations (\ref{rh4}) and (\ref{rh7}) are
thus satisfied. Furthermore, $P_{n}=O\left( h_{n}\right) $, $\Psi
_{n}=O\left( u_{n}^{-\lambda }\right) $, $\psi _{n}=O\left( h_{n}\right) $, $%
\lambda (E_{n})=O\left( 1\right) $, $\mathbb{N}_{n}=O\left( 1\right) $, $%
y_{n}=O\left( 1/\sqrt{nh_{n}^{2}}+\sqrt{h_{n}}\right) $, $L(n,\mathbf{R}%
)=O\left( h_{n}^{\gamma }\right) $, $\tau _{n}=O(h_{n}^{1/2}u_{n}^{-3\lambda
/2})$, $\alpha _{n}=O\left( u_{n}^{-3\lambda }h_{n}\log \frac{1}{h_{n}}%
\right) $, $\mathbb{M}_{n}=O\left( h_{n}\right) $. By (\ref{R1}), (\ref{R2})
and (\ref{R3}), we obtain $R_{n}(E_{n},E_{n})=O\left( h_{n}^{\gamma }\right)
$. Hence,
\[
\Omega _{n}=O\left( u_{n}^{-3\lambda }h_{n}\log \frac{1}{h_{n}}+u_{n}^{1+\mu
}+h_{n}^{\gamma }\right) .
\]
Using (\ref{rh7}), we get the rough bounds
\[
\mathbb{L}_{n}=O\left( h_{n}^{1+\gamma /d}+h_{n}\left( nh_{n}\right)
^{-(d+1)/5d}\right)
\]
and
\[
\delta _{n}=O\left( \Omega _{n}^{1/2}+h_{n}^{\gamma /d}+\left( nh_{n}\right)
^{-(d+1)/5d}+\frac{1}{nh_{n}^{2}}\right) .
\]
Note that $h_{n}^{\gamma /d}\leq h_{n}^{\gamma /2}$ and
\begin{equation}
\left( nh_{n}\right) ^{-(d+1)/5d}\leq \left( nh_{n}\right) ^{-3/10}=\left(
h_{n}^{2/5}\right) ^{3/4}\left( \left( nh_{n}^{2}\right) ^{-6/5}\right)
^{1/4}\leq h_{n}^{2/5}+\left( nh_{n}^{2}\right) ^{-6/5}.  \label{nhn}
\end{equation}
It is easy to verify that the relation (\ref{psin0}), (\ref{CC}) and (\ref
{tendzero}) are valid. Then the statement of Theorem \ref{t1.4} is valid for
\begin{eqnarray*}
0 &<&x=x_{n}=o\Big( \min \Big\{ \Big( h_{n}\log \frac{1}{h_{n}}\Big)
^{-t_{1}}l_{n}^{1/3},h_{n}^{-t_{2}}l_{n}^{1/3},\left( nh_{n}\right)
^{t_{3}}\left( \log \log \left( nh_{n}\right) \right) ^{1/3}, 
\\
&&\left. \left. h_{n}^{-\gamma /6}l_{n}^{1/3},\,\,\,\,\,\left(
nh_{n}^{2}\right) ^{1/2}\mathbf{\,}\right\} \right) ,
\end{eqnarray*}
where $t_{1}=t_{1}(\lambda ,\mu )=\frac{1+\mu }{6(1+3\lambda +\mu )}$, $%
t_{2}=t_{2}(\lambda ,\mu )=\frac{\gamma \left( 1+\mu \right) }{9\lambda }$, $%
t_{3}=t_{3}(\lambda ,\mu )=\frac{1+\mu }{21\lambda }$.

When estimating the sizes of zones of moderate deviations, we use the fact that the growth of $\varkappa _{n}$ to infinity can be arbitrarily slow. A similar remark is valid for the next examples too.

\medskip

\noindent {\bf Example 2.} Consider the density $f$ satisfying $f(x)=0$, for $x<1$%
,
\begin{equation}
C_{4}x^{-r}\leq f(x)\leq C_{5}x^{-p},\quad \mbox{for }x\geq 1,  \label{condf}
\end{equation}
and the Lipschitz condition
\[
\left\vert f(x)-f(y)\right\vert \leq C_{6}\left\vert x-y\right\vert ^{\gamma
},\quad \quad \quad \mbox{for }x\geq 1,\ y\geq 1,
\]
where $C_{4}$, $C_{5}$, $C_{6}$, $1<p\leq r<\infty $ and $0<\gamma \leq 1,$
are constants. Choose $E_{n}=\left[ 1+2h_{n},T_{n}\right] $, where
\[
\begin{array}{c}
T_{n}=\min \Big\{ \left( h_{n}\log \frac{1}{h_{n}}\right) ^{-\frac{1}{3r+p-1%
}},h_{n}^{-\gamma /p},\left( nh_{n}^{2}\right) ^{3/(p+5)}, \\
\left. \varkappa _{n}\,h_{n}^{-2\gamma /3r},\varkappa _{n}\,\left(
nh_{n}\right) ^{2/7r}\right\} ,
\end{array}
\]
and $\varkappa _{n}$ is an arbitrary sequence of positive numbers which tend
to zero as $n\rightarrow \infty $ slower than a fixed function of $h_{n}$.
Without loss of generality we assume $h_{n}\leq 1/8$. Then it is easy to
estimate (sometimes roughly) $\phi _{n}=O\left( T_{n}^{1-p}\right) $, $\beta
_{n}^{-1}=O\left( T_{n}^{r}\right) $, $q_{n}(x)=O\left( h_{n}^{\gamma
}\right) $, $D_{n}=O\left( 1\right) $, $\varepsilon _{n}=O\left(
h_{n}^{\gamma }\right) $, $Q_{n}=O\left( h_{n}^{\gamma }T_{n}^{r}\right)
=o\big( h_{n}^{\gamma /3}\big) $. The relations (\ref{rh4}) and (\ref{rh7}%
) are thus satisfied. Furthermore, $P_{n}=O\left( h_{n}\right) $, $\Psi
_{n}=O\left( T_{n}^{r}\right) $, $\psi _{n}=O\left( h_{n}\right) $, $\lambda
(E_{n})\leq T_{n}$, $\mathbb{N}_{n}=O\left( 1\right) $,
\[
y_{n}=O\bigg( \frac{T_{n}}{\left( nh_{n}^{2}\right) ^{1/2}}+\sqrt{h_{n}}%
\bigg) =O\left( \left( nh_{n}^{2}\right) ^{\left( 1-p\right) /2(p+5)}+\sqrt{%
h_{n}}\right) ,
\]
\[
L(n,\mathbf{R})=O\left( \int_{E_{n}}q_{n}(x)\,dx\right) =O\left(
h_{n}^{\gamma }T_{n}\right) ,
\]
\[
\tau _{n}=O(T_{n}^{3r/2}h_{n}^{1/2}),\quad \alpha _{n}=O\left(
T_{n}^{3r}h_{n}\log \frac{1}{h_{n}}\right) ,
\]
\[
R_{n}(E_{n},E_{n})=O\left( h_{n}+\int_{E_{n}}q_{n}(x)\,dx\right) =O\left(
h_{n}^{\gamma }T_{n}\right) ,
\]
\[
\Omega _{n}=O\left( T_{n}^{3r}h_{n}\log \frac{1}{h_{n}}+T_{n}^{1-p}+h_{n}^{%
\gamma }T_{n}\right) .
\]
To estimate $\mathbb{L}_{n}$ in (\ref{Ln}), we begin with estimating (for $n$
large enough) the integral
\begin{eqnarray}
\int_{E_{n}}\Delta _{n}^{2/5}(x)\,\left( q_{n}(x)/f(x)\right) ^{1/d}f(x)\,dx
&\leq &\int_{E_{n}}\frac{h_{n}^{\gamma /d}}{\left( \beta _{n}\,nh_{n}\right)
^{1/5}\beta _{n}^{1/2}}\,\left( f(x)\right) ^{1/2-1/d+1}\,dx  \nonumber \\
&\leq &h_{n}^{\gamma /d}\int_{E_{n}}\left( f(x)\right) ^{1/2}\,dx\leq
h_{n}^{\gamma /d}T_{n}^{1/2}  \label{Ln1}
\end{eqnarray}
(using $1\leq d\leq 2$, (\ref{tendzero}), and (\ref{condf})). Moreover,
\begin{eqnarray}
&&\int_{E_{n}}\Delta _{n}^{2(d+1)/5d}(x)\,f(x)\,dx  \nonumber \\
&\leq &\left\{
\begin{tabular}{l}
$O\left( \left( nh_{n}\right) ^{-(d+1)/5d}\right) ,\quad $ $\,$ $\mbox{for }%
(4d-1)p>5,$ \\
$O\left( \left( nh_{n}\right) ^{-(d+1)/5d}\log \frac{1}{h_{n}}\right) ,\quad
$ $\,$ $\mbox{for }(4d-1)p=5,$ \\
$O\left( \left( nh_{n}\right) ^{-(d+1)/5d}T_{n}^{1-(4d-1)p/5d}\right) ,\quad
$ $\,$ $\mbox{for
}(4d-1)p<5.$%
\end{tabular}
\ \ \ \ \right.  \label{Ln22}
\end{eqnarray}
By (\ref{Ln}), (\ref{Ln1}) and (\ref{Ln22}), we have
\begin{equation}
\mathbb{L}_{n}=\left\{
\begin{tabular}{l}
$O\left( h_{n}^{1+\gamma /d}T_{n}^{1/2}+h_{n}\left( nh_{n}\right)
^{-(d+1)/5d}\right) ,\quad $ $\,$ $\mbox{for }(4d-1)p>5,$ \\
$O\left( h_{n}^{1+\gamma /d}T_{n}^{1/2}+h_{n}\left( nh_{n}\right)
^{-(d+1)/5d}\log \frac{1}{h_{n}}\right) ,\quad $ $\,$ $\mbox{for }(4d-1)p=5,$
\\
$O\left( h_{n}^{1+\gamma /d}T_{n}^{1/2}+h_{n}\left( nh_{n}\right)
^{-(d+1)/5d}T_{n}^{1-(4d-1)p/5d}\right) ,\quad $ $\,$ $\mbox{for
}(4d-1)p<5,$%
\end{tabular}
\ \ \ \ \right.  \label{lnis}
\end{equation}
and $\mathbb{M}_{n}=O\left( h_{n}T_{n}\right) $. Applying (\ref{lnis}) with $%
d=2$ and using (\ref{nhn}), we obtain
\[
\delta _{n}=O\left( \Omega _{n}^{1/2}+h_{n}^{\gamma
}T_{n}+h_{n}^{2/5}+\left( nh_{n}^{2}\right) ^{\left( 1-p\right)
/(p+5)}\right) .
\]
The relations (\ref{psin0}), (\ref{CC}) and (\ref{tendzero}) are again
valid. Then the statement of Theorem \ref{t1.4} is valid for
\begin{multline}
0 <x=x_{n}=o\Big( \min \Big\{ \Big( h_{n}\log \frac{1}{h_{n}}\Big)
^{-t_{4}}l_{n}^{1/3},h_{n}^{-t_{5}}l_{n}^{1/3},h_{n}^{-t_{6}}l_{n}^{1/3},%
 \\
\left. \left. \,\,\left( nh_{n}^{2}\right) ^{\left( p-1\right)
/2(p+5)}\right\} \right) ,
\end{multline}
where $t_{4}=t_{4}(p,r)=\frac{p-1}{6(3r+p-1)}$, $t_{5}=t_{5}(p,r,\gamma )=%
\frac{\gamma \left( p-1\right) }{9r}$, $t_{6}=t_{6}(p,r,\gamma )=\frac{%
\gamma \left( p-1\right) }{6p}$.

\medskip

\noindent {\bf Example 3.} Consider the density $f$ satisfying $f(x)=0$
for $x\notin \lbrack 0,1]$,
\[
C_{7}\leq f(x)\leq C_{8}x^{-v}\quad \mbox{for }0<x\leq 1,
\]
and
\[
\left\vert \mathbf{\,}f(x)-f(y)\right\vert \leq C_{9}\left( \left\vert
x-y\right\vert ^{\gamma }x^{-v}+\left\vert x-y\right\vert \,x^{-1-v}\right)
,\quad \mbox{for }0<x\leq y\leq 1,
\]
where $C_{7}$, $C_{8}$, $C_{9}$, $0<\gamma <1$, $0<v<1$ are constants. For
example, such a condition is satisfied if
\[
f(x)=x^{-v}g(x),\quad \mbox{for }0<x\leq 1,
\]
where the function $g$ satisfies the H\"{o}lder condition with parameter $%
\gamma $.\ Choose $E_{n}=\left[ w_{n},1-2h_{n}\right] $, where
\[
w_{n}=\max \Big\{ \Big( h_{n}^{1-v}\log \frac{1}{h_{n}}\Big) ^{\frac{1}{%
2v+1}},h_{n}^{3\gamma
/(1+2v)}l_{n}^{2/(1+2v)},h_{n}^{3/(2v+4)}l_{n}^{1/(v+2)},\varkappa
_{n}\,h_{n}^{2\gamma /3v}\Big\} ,
\]
and $\varkappa _{n}\rightarrow 0$ as $n\rightarrow \infty $ is an arbitrary
sequence of positive numbers which tend to infinity slower than a fixed
function of $h_{n}^{-1}$. Without loss of generality we assume $h_{n}\leq
1/8 $. Then we estimate $\phi _{n}=O\left( w_{n}^{1-v}\right) $, $\beta
_{n}^{-1}=O\left( 1\right) $, $D_{n}=O\left( w_{n}^{-v}\right) $, $%
q_{n}(x)=O\left( h_{n}^{\gamma }\,x^{-v}+h_{n}\,x^{-v-1}\right) $, for $%
3h_{n}/2<x<1-3h_{n}/2$, and $Q_{n}=O\left( h_{n}^{\gamma
}w_{n}^{-v}+h_{n}w_{n}^{-v-1}\right) $. The relations (\ref{rh4}) and (\ref
{rh7}) are thus satisfied. By (\ref{SV}), we have $\varepsilon _{n}=O\left(
h_{n}^{\gamma }w_{n}^{-v}+h_{n}w_{n}^{-v-1}\right) $. Furthermore, $%
P_{n}=O\left( h_{n}^{1-v}\right) $, $\Psi _{n}=O\left( w_{n}^{-v}\right) $, $%
\psi _{n}=O\left( h_{n}w_{n}^{-v}\right) $ and
\[
\mathbb{N}_{n}=O\left( \mathbf{1}\left\{ v<2/3\right\} +\mathbf{1}\left\{
v=2/3\right\} \log \frac{1}{h_{n}}+\mathbf{1}\left\{ v>2/3\right\} \mathbf{\,%
}w_{n}^{1-3v/2}\right) .
\]
To estimate $w_{n}^{1-3v/2}$, we note that $w_{n}\geq h_{n}$ and, hence, $%
w_{n}^{1-3v/2}\leq h_{n}^{-1/2}w_{n}^{1-v}$. Thus,
\[
y_{n}=O\left( 1/\sqrt{nh_{n}^{2}}+\mathbb{N}_{n}\sqrt{h_{n}}\right) =O\left(
1/\sqrt{nh_{n}^{2}}+w_{n}^{1-v}\right) .
\]
To obtain a rough bound for $\mathbb{L}_{n}$ in (\ref{Ln}), we begin with
estimating (for $n$ large enough) the integral
\begin{eqnarray}
&&\int_{E_{n}}\Delta _{n}^{2/5}(x)\,\left( q_{n}(x)/f(x)\right)
^{1/d}f(x)\,dx  \nonumber \\
&\leq &\frac{\left( O\left( h_{n}^{\gamma
}w_{n}^{-v}+h_{n}w_{n}^{-v-1}\right) \right) ^{1/d}}{\left( nh_{n}\right)
^{1/5}}\,  \nonumber \\
&\leq &O\left( h_{n}^{\gamma }w_{n}^{-v}+h_{n}w_{n}^{-v-1}\right) +\frac{1}{%
\left( nh_{n}\right) ^{2/5}}.  \label{Ln7}
\end{eqnarray}
Moreover,
\begin{equation}
\int_{E_{n}}\Delta _{n}^{2(d+1)/5d}(x)\,f(x)\,dx=O\left( \left(
nh_{n}\right) ^{-(d+1)/5d}\right) .  \label{Ln8}
\end{equation}
Using (\ref{Ln}), (\ref{Ln7}) and (\ref{Ln8}), we see that
\begin{equation}
\mathbb{L}_{n}=O\left( h_{n}\left( h_{n}^{\gamma
}w_{n}^{-v}+h_{n}w_{n}^{-v-1}+\left( nh_{n}\right) ^{-(d+1)/5d}\right)
\right) .  \label{Ln9}
\end{equation}
By (\ref{rh4}), $\mathbb{M}_{n}=O\left( h_{n}\right) $. Furthermore, $%
\lambda (E_{n})=O\left( 1\right) $, $R_{n}(E_{n},E_{n})=O(h_{n}^{\gamma
}+h_{n}w_{n}^{-v})$, $\tau _{n}=O(w_{n}^{-3v/2}h_{n}^{(1-v)/2})$, $\alpha
_{n}=O\left( w_{n}^{-3v}h_{n}^{1-v}\log \frac{1}{h_{n}}\right) $. By (\ref
{SV}), we have $L(n,\mathbf{R})=O(h_{n}^{\gamma }+h_{n}w_{n}^{-v})$.
Collecting our bounds and using (\ref{nhn}), we obtain $\delta _{n}=O\left(
\Omega _{n}^{1/2}+h_{n}^{\gamma }w_{n}^{-v}+h_{n}w_{n}^{-v-1}+h_{n}^{2/5}+%
\frac{1}{nh_{n}^{2}}\right) $, $$\Omega _{n}=O\left(
w_{n}^{-3v}h_{n}^{1-v}\log \frac{1}{h_{n}}+w_{n}^{1-v}+h_{n}^{\gamma
}\right) .$$ Note that
\[
\frac{w_{n}^{-3v}}{nh_{n}}=\frac{h_{n}^{1-v}}{w_{n}^{1+2v}}%
w_{n}^{1-v}h_{n}^{v}\frac{1}{nh_{n}^{2}}=o(1)
\]
and, hence, $D_{n}^{1/2}/\left( nh_{n}\right) ^{1/5}=o(1)$. Now it is easy
to verify that the relations (\ref{psin0}), (\ref{CC}) and (\ref{tendzero})
are valid.

The statement of Theorem \ref{t1.4} is valid for
\begin{eqnarray*}
0 &<&x=x_{n}=o\Big( \min \Big\{ \Big( h_{n}^{1-v}\log \frac{1}{h_{n}}%
\Big)
^{-t_{7}}l_{n}^{1/3},h_{n}^{-t_{8}}l_{n}^{t_{9}},h_{n}^{-t_{10}}l_{n}^{t_{11}},
\\
&&\left. \left. h_{n}^{-t_{12}}l_{n}^{1/3},\,\,\,\,\left( nh_{n}^{2}\right)
^{1/2}\right\} \right) ,
\end{eqnarray*}
where $t_{7}=\frac{1-v}{6(1+2v)}$, $t_{8}=\frac{\gamma (1-v)}{2(1+2v)}$, $%
t_{9}=\frac{v}{1+2v}$, $t_{10}=\frac{1-v}{4(v+2)}$, $t_{11}=\frac{v+1}{2(v+2)%
}$, $t_{12}=\frac{\gamma (1-v)}{9}$. $\medskip $

We shall say that the density $g$ belongs to the class $\mathcal{D}_{k}$ if
it may be represented in the form
\begin{equation}
g(x)=f(b(x-a)),  \label{gf}
\end{equation}
where $a$ and $b\neq 0$ are real constants and $f$ is a density satisfying
the conditions of Example k. It is clear that if we replace the density $f$
by the density $g$ from (\ref{gf}) in Example k, then the formula for $%
x=x_{n}$ from the corresponding statement of Theorem \ref{t1.4} will be the
same.\medskip

\noindent {\bf Example 4.} Consider the density $f$ of the form $%
f(x)=\sum_{j=1}^{m}p_{j}f_{j}(x)$, where $0\leq p_{j}\leq 1$, $%
\sum_{j=1}^{m}p_{j}=1$, the functions $f_{j}=r_{j}\,1_{\mathcal{J}_{j}}$ are
densities belonging to $\cup _{k=1}^{3}\mathcal{D}_{k}$ and $%
J_{j}=(x_{j},y_{j})$, $x_{j}<y_{j}$, $j=1,2,\ldots ,m$, is a finite
collection of disjoint intervals. It is obvious that we should choose the
sets $E_{n}$ as $E_{n}=\bigcup_{j=1}^{m}E_{n}^{(j)}$, where $E_{n}^{(j)}$ is
the set $E_{n}$ which corresponds to the density $f_{j}$. Then the statement
of Theorem \ref{t1.4} is valid for
\[
0<x=x_{n}=o\Big( \min_{1\leq j\leq m}x_{n}^{(j)}\Big) ,
\]
where $x_{n}^{(j)}$ is $x_{n}$ which corresponds to the density $f_{j}$%
.

\end{document}